\documentclass[11pt]{article} 
\usepackage{amssymb}
\usepackage{amsmath}
\usepackage{amsthm}
\newtheorem{prop}{\bfseries Proposition}
\newtheorem{lem}{\bfseries Lemma}
\newtheorem{thm}{\bfseries Theorem}
\newtheorem{cor}{\bfseries Corollary}
\newcommand{\dil}[2]{\left\langle #1 \right\rangle_{\! #2}} 
\def\SX{{\mathsf X}}
\def\SY{{\mathsf Y}}
\def\SZ{{\mathsf Z}}
\def\ST{{\mathsf T}}
\def\SF{{\mathsf F}}
\def\SR{{\mathsf R}}
\def\SM{{\mathsf M}}
\def\CR{{\mathcal R}}
\def\CF{{\mathcal F}}
\begin{document}

\begin{center}
{\bf\Large\sc 
Tetrahedron equation, Weyl group,\\
 and quantum dilogarithm  
} 

\vspace*{4mm}
{\sc Andrei Bytsko \ and \ Alexander Volkov }
 
\vspace*{5mm}
{\em Dedicated to Ludwig Faddeev 
on the occasion of his 80th birthday }

\end{center}
\vspace*{2mm}

\begin{abstract}
We derive a family of solutions to
the tetrahedron equation using the RTT  
presentation of a two parametric 
quantized algebra of regular functions on an 
upper triangular subgroup of $GL(n)$.
The key ingredients of the construction
are the longest element of the Weyl group,
the quantum dilogarithm function, and 
central elements of the quantized 
division algebra of rational functions
on the subgroup in question.

\end{abstract}

\section{Introduction}

In what follows, 
$q,\lambda,\mu,\nu \in {\mathbb C}$
and $|q|\neq 1$. All algebras are considered
over~$\mathbb C$.

The {\em quantum exponential function}
is the following formal power series in $x$:
\begin{equation}\label{Sx}
  \dil{x}{q} =  \sum_{n = 0}^{\infty}  
   \frac{(-x)^n }{(1-q)\ldots(1-q^{n})} \,.
\end{equation} 
$\dil{x}{q}$ is the
unique series in ${\mathbb C}[[x]]$
which satisfies the functional relation
\begin{equation}\label{fundil}  
  \dil{q\,x}{q} = (1+x) \, \dil{x}{q} 
\end{equation}
and the condition $\dil{0}{q}=1$.
The series inverse to
$\dil{x}{q}$ is given by
\begin{equation}\label{Sxi}
  1/\dil{x}{q} = \dil{q^{-1} x}{q^{-1}} =
   \sum_{n = 0}^{\infty}  
   \frac{ q^{\frac{n(n-1)}{2}} x^n }{(1-q)\ldots(1-q^{n})} \,.
\end{equation} 
Indeed, using (\ref{fundil}), one verifies that
$\theta(x) \equiv \dil{x}{q} \dil{q^{-1} x}{q^{-1}}$
has the property $\theta(q x)=\theta(x)$.

Quantum exponential function satisfies also 
the following identities \cite{Sch1,FK1,FV1}:
\begin{align}
\label{qexp} 
{}&  \dil{\SX}{q} \,  \dil{\SY}{q}  =  
 \dil{ \SX + \SY}{q} \,, \qquad\quad 
 \dil{\SX}{q} \, \dil{\SX\SY}{q} \, \dil{\SY}{q}  =  
 \dil{\SY}{q} \, \dil{\SX}{q} \,, 
\end{align}
involving two $q$--commuting indeterminates,
$\SY \SX =q\, \SX \SY$. The second identity 
in (\ref{qexp}) can be regarded  
as a quantum counterpart of the five term
relation for the Rogers dilogarithm, see~\cite{FK1}.
For this reason, $\dil{x}{q}$ was named the
{\em quantum dilogarithm}.

Using the quantum dilogarithm identity twice, one can
derive another identity \cite{KV2}:
\begin{align}
\label{qexp4} 
{}&   
 \dil{\SX}{q} \, \dil{\SX\SY}{q} \, \dil{\SZ}{q} \,
 \dil{\SY}{q}  =  \dil{\SZ}{q} \, \dil{\SZ\SX}{q} \,
 \dil{\SY}{q} \, \dil{\SX}{q} \,,  
\end{align}
involving three pairwise $q$--commuting indeterminates,
\begin{equation}\label{XYZ}
\SY \SX = q \, \SX \SY \,, \qquad
\SX \SZ = q \, \SZ \SX \,, \qquad
\SZ \SY = q \, \SY \SZ \,.
\end{equation} 

Below we will need a more general statement. Namely,
\begin{lem}\label{pentdil}
Let indeterminates $\SX$, $\SY$, $\SZ$   
satisfy relations~(\ref{XYZ}). 
Let $S_{q,\lambda}(t), U_{q,\lambda}(t)
\in {\mathbb C}[[t]]$ 
be non--constant formal power series in $t$ with
coefficients depending on $q$ and~$\lambda$. 
Then the following relations hold
\begin{align}
\label{pent}
{} S_{q,\lambda}(\SX) \, S_{q,\lambda\mu}(\SX\SY) \,  
 S_{q,\mu}(\SY) &=
 S_{q,\mu}(\SY) \, S_{q,\lambda}(\SX) \,, \\[1mm]
{} 
\label{tetr}
 U_{q,\lambda}(\SX) \, U_{q,\lambda\mu}(\SX\SY) \, 
 U_{q,\nu}(\SZ) \,  U_{q,\mu}(\SY) &=
 U_{q,\nu}(\SZ) \, U_{q,\lambda\nu}(\SZ\SX) \, 
 U_{q,\mu}(\SY) \, U_{q,\lambda}(\SX)  
\end{align}
if and only if $S_{q,\lambda}(t)$ and $U_{q,\lambda}(t)$ 
 are given by
\begin{align}
\label{diln}
{}& 
 S_{q,\lambda}(t) = 
 \dil{q^{\frac{m(m-1)}{2}} \lambda^k \, t^m}{q^{m^2}}  \,,\qquad 
 U_{q,\lambda}(t) = \kappa(q)\, 
 \dil{q^{\frac{m(m-1)}{2}} \lambda^k \, t^m}{q^{m^2}}  \,,
\end{align}
where $m \in {\mathbb Z}_+$, $k \in {\mathbb Z}$,
and $\kappa(q)$ is an arbitrary function of~$q$.
\end{lem}

Note that the order of the middle factors in (\ref{qexp4}) 
and (\ref{tetr}) can be reversed because 
$\SX\SY$ commutes with $\SZ$ and likewise
$\SZ\SX$ commutes with~$\SY$. Let us also remark
that the commutative parameters $\lambda$, $\mu$, $\nu$
appearing in (\ref{pent})--(\ref{tetr}) can be 
set equal to unity by a suitable rescaling of  
 $\SX$, $\SY$,~$\SZ$.
Therefore, these parameters are not of the same
nature as the ``physical'' spectral parameters 
that appear, for instance, in the tetrahedron
equation in~\cite{Za}. However, we   
keep these parameters because equations which we 
consider remain nontrivial when some of these
parameters are sent to zero. For instance,
if $k$ in (\ref{diln}) is positive,
equation (\ref{pent}) can be obtained from 
(\ref{tetr}) by setting $\nu=0$ and taking 
into account that $\dil{0}{q}=1$.

The quantum exponential function often appears
as a building block for solutions to the
{\em Yang--Baxter equation},
\begin{equation}\label{RYB}
 R_{12} R_{13} R_{23}  =
 R_{23} R_{13} R_{12}   \,.
\end{equation}
See, for instance, Drinfeld's solution \cite{Dri}
for the R--matrix of $U_q(sl_2)$. 

The Yang--Baxter equation has an interpretation
in terms of three pairwise crossing straight lines in a plane.
Its higher dimensional counterpart involving
four pairwise intersecting planes in ${\mathbb R}^3$ is 
Zamolodchikov's {\em tetrahedron equation} \cite{Za}:
\begin{equation}\label{Rtetr}
 R_{123} R_{145} R_{246} R_{356} =
 R_{356} R_{246} R_{145}  R_{123} \,.
\end{equation}

The symbol $R$ in (\ref{RYB}) and (\ref{Rtetr}) can be 
assigned different meanings. 
For instance, it can be a set--theoretic map (and 
then the product of $R$'s is understood as a composition) or 
a linear operator on a vector space.
In the present article, we will regard $R$ as a function
of non--commuting indeterminates.

In \cite{KV1} (see also \cite{KS1} for more details
and \cite{KO} for some generalizations),
a solution to the tetrahedron equation
 was constructed using the representation
theory of $GL_q(n)$ --- the quantized algebra of regular functions
on the group $GL(n)$. Let us recall the scheme of this
solution. Let $\pi$ be the irreducible representation
of $GL_q(2)$ on an infinite--dimensional space $V$ 
constructed in~\cite{SoV1}. The algebra 
homomorphisms $\varphi_k : GL_q(n) \to GL_q(2)$, where
$k \in [1,n-1]$ labels a vertex of the $A_{n-1}$ Dynkin diagram,
correspond to the embeddings 
$GL(2) \subset GL(n)$ of the classical groups.
Let $s_{i_1}\cdots s_{i_k}$ be a reduced form of the 
longest element of the Weyl group for $GL(n)$.
Then, as was shown in \cite{SoV2}, 
$\pi_{i_1,\ldots,i_k} : x \to \bigl((\pi \circ \varphi_{i_1})\otimes
 \cdots \otimes (\pi \circ \varphi_{i_k})\bigr) \circ 
 (\Delta \otimes \cdots \otimes \Delta)(x)$
is an irreducible representation of $GL_q(n)$.
Here $\Delta$ is the comultiplication of $GL_q(n)$.
Let $S \in \text{End}\, V^{\otimes 3}$ be the intertwiner 
of two irreducible representations
of $GL_q(3)$: $S\, \pi_{121}(x) = \pi_{212}(x)\, S$,
and let $P_{13}$ denote the permutation
of the first and the third tensor factors in $V^{\otimes 3}$.
It was observed in \cite{KV1} that $R = S \circ P_{13}$
is a solution to the tetrahedron equation 
because both sides of (\ref{Rtetr}) for
such $R$ are intertwiners of the representations
$\pi_{121321}$ and $\pi_{323123}$ of $GL_q(4)$.
Since these representations are irreducible, the
intertwiners must coincide if they coincide on
some vector from $V^{\otimes 6}$.

In the outlined above construction, 
the l.h.s. and 
the r.h.s. of (\ref{Rtetr}) correspond
to two ways of transforming the longest element
of the  Weyl group for $A_3$ 
with the help of the braid relations 
$s_{k} s_{k'} s_k = s_{k'} s_{k} s_{k'}$
if $|k-k'|=1$ and 
$s_k s_{k'} = s_{k'} s_{k}$ if $|k-k'|>1$:
\begin{equation}\label{We123}
\begin{aligned}
{}& s_1 s_2 s_1 s_3 s_2 s_1 =  
 s_2 s_1 s_2 s_3 s_2 s_1 =
 s_2 s_1 s_3 s_2 s_3 s_1 =
 s_2 s_3 s_2 s_1 s_2 s_3 =
 s_3 s_2 s_1 s_3 s_2 s_3  \,, \\
{}& s_1 s_2 s_1 s_3 s_2 s_1 = 
  s_1 s_2 s_3 s_2 s_1 s_2 =
  s_1 s_3 s_2 s_3 s_1 s_2 =
  s_3 s_2 s_1 s_2 s_3 s_2  = 
  s_3 s_2 s_1 s_3 s_2 s_3 \,.
\end{aligned} 
\end{equation} 

The aim of the present article is to construct
solutions to (\ref{Rtetr}) based on the same Weyl group 
structure but without utilising irreducible representations.
For this purpose, we will first search for $R$ that
transforms a $s_1 s_2 s_1$ product of non--commuting matrices 
into a $s_2 s_1 s_2$ product. 
This step is similar to finding the so--called fundamental 
solution to the
Yang--Baxter equation by solving the RLL relation.
An important difference of our approach from the  
above--mentioned construction using irreducible representations 
of $GL_q(n)$ is that 
we will deal with a division algebra  
whose center is not scalar.
As a consequence, $R$ is not determined uniquely by the
intertwining relation. Our first principal result
is that the part of $R$ which is fixed by the 
intertwining relation is 
the quantum dilogarithm function of 
the image $w$ of a central element of the division algebra 
whereas the ambiguity factor in~$R$ is a function 
of the dual image $\tilde{w}$ of the same central element
(it turns out that $\tilde{w}$ and $w$ commute),
cf.~Theorem~\ref{RFG}. 
Further, motivated by the obtained expression for $R$,
we will search for a solution to (\ref{Rtetr}) of the form
$\tilde{R}=\SF \,f(x)\,g(y)$, where $x$ and $y$  are monomials
in the generators of a $q$--Weyl algebra and
$\SF$ is a certain involution (a counterpart of the permutation
in the Yang--Baxter equation) such that $\SF\,w=\tilde{w}\,\SF$. 
Our second principal result is that both functions $f$ and $g$ 
are required to
be the quantum dilogarithm function and $x$, $y$ must
coincide with $w$ and $\tilde{w}$, cf.~Theorem~\ref{RS}. 
This implies, in particular, that $R$ found by 
solving the intertwining relation satisfies the tetrahedron 
equation if the ambiguity factor is also the quantum 
dilogarithm function.
We observe that the latter condition holds automatically 
if $R$ is required to be an involution. 

The proofs of all statements are given in the Appendix.

\section{An upper triangular quantum group}

In order to construct a family of solutions to the
tetrahedron equation
we will use (a certain quotient of) the two--parameter quantum 
group $GL_{q,p}(n)$.
The following two--parameter generalization
of the Jimbo--Drinfeld R--matrix is well known~\cite{Re1}:
\begin{equation}\label{Rpq}
 R_{q,p}^{(n)} = 
 qp \sum_{i=1}^n E_{ii} \otimes E_{ii} +
 q\!  \sum_{1 \leq i <  j \leq n} E_{ii} \otimes E_{jj} +
 p\!  \sum_{1 \leq j < i \leq n} E_{ii} \otimes E_{jj} +
 (qp -1)\! \sum_{1 \leq i <  j \leq n} E_{ij} \otimes E_{ji} ,
\end{equation}
where $E_{ij}$ are the basis matrices, i.e. 
$(E_{ij})_{kl}=\delta_{ik} \delta_{jl}$. 
Matrix $R_{q,p}^{(n)}$ satisfies the Yang--Baxter equation~(\ref{RYB}).
Matrix ${\cal P}^{(n)} = \sum_{i,j=1}^n E_{ij} \otimes E_{ji}$ 
acts on ${\mathbb C}^n \otimes {\mathbb C}^n$ as
the permutation of the tensor factors 
and also provides a solution to~(\ref{RYB}).

Using the RTT approach \cite{FRT}, 
we define an associative algebra 
$GL^+_{q,p}(n)$  
by the following matrix presentation:
\begin{equation}\label{RGG}
 \check{R}_{q,p}^{(n)} \,  \bigl(X(n) \otimes X(n) \bigr) =
 \bigl(X(n) \otimes X(n) \bigr) \, \check{R}_{q,p}^{(n)} \,,
\end{equation} 
where 
$\check{R}_{q,p}^{(n)} \equiv {\cal P}^{(n)} R_{q,p}^{(n)}$,
$\otimes$ stands for the Kronecker product, and 
 $X(n)$ is an $n\,{\times}\,n$ upper triangular 
matrix with non--commuting entries,
$x_{ij}$, $1 \leq i \leq j \leq n$.

Explicitly, for  
every quadruple 
($x_{ij}$, $x_{ik}$, $x_{mj}$, $x_{mk}$) where $i<m \leq j<k$, 
eq. (\ref{RGG}) is
equivalent to the defining relations of $GL_{q,p}(2)$ :
\begin{equation}\label{G2a}
\begin{aligned}
{}&  x_{ij}x_{ik} =q\, x_{ik}x_{ij} , \quad 
 x_{mj}x_{mk} =q\, x_{mk}x_{mj} , \quad
  x_{ij}x_{mj} = p\, x_{mj}x_{ij} , \quad 
  x_{ik}x_{mk} = p\, x_{mk}x_{ik} , 
  \\  
{}&  x_{ik}x_{mj} = \frac{p}{q}\, x_{mj}x_{ik} \,, \qquad
 x_{ij}x_{mk} - x_{mk}x_{ij} = (q-\frac{1}{p})\, x_{ik}x_{mj} =
  (p-\frac{1}{q})\, x_{mj}x_{ik} \, .
\end{aligned}
\end{equation} 
For every triple
($x_{ij}$, $x_{ik}$, $x_{mk}$), where
$i<m$, $j<k$, $j<m$, eq. (\ref{RGG}) is
equivalent to
\begin{align}\label{B2}
{}&  x_{ij}\,x_{ik} =q\, x_{ik}\,x_{ij} \,, \qquad  
 x_{ik}\,x_{mk} = p\, x_{mk}\,x_{ik} \,, 
  \qquad  x_{ij}\,x_{mk} = x_{mk}\,x_{ij}   \,,
\end{align}
which is a reduction of (\ref{G2a}) 
obtained by setting $x_{mj}$ to zero.

In what follows, we will consider a special case where
$p$ is an integer power of~$q$, 
\begin{align}\label{pqq}
  p = q^{r} \,, \qquad  
   r \in {\mathbb Z} \,, \qquad r \neq -1 \,.
\end{align}

\begin{lem}\label{Cdet}
If $p$ and $q$ are related as in (\ref{pqq}), then
the element
\begin{align}\label{c0}
 C_0 = \prod_{i=1}^n (x_{ii})^{r^{n-i}} 
\end{align}  
is central in $GL_{q,p}^+(n)$.
\end{lem}

Let $SL_{q,q^r}^+(n)$ denote the quotient of 
$GL^+_{q,q^r}(n)$ by the ideal 
generated by $(C_0 -1)$.
We use the $SL$ notation
for the resulting algebra since, in the $p=q$ 
case, $C_0$ coincides with the determinant of $X(n)$.
Let us remark that $SL_{q,q^r}^+(n)$
can be equipped with the structure of a bialgebra
by defining the comultiplication in the standard
way, i.e.,
$\Delta \bigl(X(n)\bigr) =  X(n) \otimes X(n)$.
The element $C_0$ is group like w.r.t. this
comultiplication.

Along with the algebra of regular functions
$SL_{q,q^r}^+(n)$ generated by $x_{ij}$
we will consider the
corresponding division algebra $D_{q,r}(n)$
of rational functions in non--commuting 
indeterminates~$x_{ij}$ satisfying relations
(\ref{G2a})--(\ref{B2}).
Note that, in $D_{q,r}(n)$,
the constraint $C_0=1$ can be
resolved by expressing $x_{nn}$ in terms of other~$x_{ii}$
which decreases the number of independent
indeterminates by one.
In particular, for $r=0$, we have $x_{nn}=C_0=1$.

In general, for $D_{q,r}(2)$,
we have $x_{22}=x_{11}^{-r}$ and the only
defining relations is
\begin{align}\label{B22}
{}&  x_{11} \, x_{12} =q\, x_{12}\,x_{11}  \,.
\end{align} 

For $D_{q,r}(3)$, we have 
$x_{33} = x_{11}^{-r^2} x_{22}^{-r} $ 
along with the following set of 
relations for the independent indeterminates: 
\begin{align}\label{B3}
{}&  x_{11} x_{12}=q\, x_{12} x_{11}  , \quad  
 x_{11} x_{13}=q\, x_{13} x_{11} , \quad  
 x_{12} x_{13}=q\, x_{13} x_{12} ,  \quad
 x_{22} x_{23}=q\, x_{23} x_{22} ,
  \\
{}& x_{13} x_{23} =  q^r x_{23} x_{13}  , \quad
 x_{12} x_{22} =  q^r x_{22} x_{12}  ,  \quad
 x_{13} x_{22} =  q^{r-1} x_{22} x_{13} , \\
\label{B33}
{}&    
 x_{12} x_{23} - x_{23} x_{12} =
 (q -q^{-r})\, x_{13} x_{22}  . 
\end{align}  

A remark is in order here.
Let $SL^+_r(n)$ denote the subgroup
of $n\,{\times}\,n$ upper triangular matrices with
commuting entries satisfying the condition $C_0=1$,
where $C_0$ is given by~(\ref{c0}). 
The quasi--classical limit,
$q=e^{\hbar}$, $p=e^{r\hbar}$, $\hbar \to 0$,
of the RTT relations (\ref{RGG}) for $X(n)$ induces the structure
of a Poisson--Lie group on $SL^+_r(n)$.
In particular, the classical counterpart of 
relation (\ref{B22}) is 
a log--canonical Poisson bracket:
$\{x_{11},x_{12}\} =  x_{11} x_{12}$.
For $n=3$, the Poisson structure corresponding
to (\ref{B3})--(\ref{B33}) is necessarily 
degenerate since it involves five variables
whereas the dimension
of the maximal symplectic leaf must be even and thus
it is four. Therefore, the algebra $D_{q,r}(3)$
must have a non--trivial central element. 
 
For $GL_q(n)$, the center of the
corresponding division algebra was described 
in~\cite{Cli,Pa1}. In the multi--parameter case,
the structure of the center depends in general on the 
values of parameters~\cite{MoPa}. However, for $D_{q,r}(3)$,
it is not hard to find the center by a direct computation.
Indeed, using (\ref{B3})--(\ref{B33}), it is easy to check 
that the subalgebra of
elements commuting with $x_{11}$, $x_{22}$, 
and $x_{13}$ is generated by 
$x^{r-1}_{11} x_{12} x_{23} x_{13}^{-1}$
and $x^{r-1}_{11} x_{22}$.
Imposing the condition of commutativity
with $x_{12}$ and $x_{23}$, we deduce that
the center of $D_{q,r}(3)$ is generated by the element
\begin{align}\label{Dpq2} 
  C_1 =  (x_{12} x_{23} - q x_{13} x_{22})\, 
   x_{11}^{-1}  x_{33}^{-1} x_{13}^{-1} 
   = 
   (x_{12} x_{23} - q x_{13} x_{22})\, 
   x_{11}^{r^2-1} x_{22}^r x_{13}^{-1}  \, .
\end{align}

\section{Main construction}

Let $Q_q(n)$ denote the $q$--Weyl algebra 
with generators $a_i, a_i^{-1}, b_i, b_i^{-1}$, $i=1,{\ldots},n$
and defining relations 
\begin{equation}\label{B222}
\begin{aligned}
{}&  a_i \,b_j  =q^{\delta_{ij}} \, b_j \,a_i \,, \qquad
 a_i \,a_j =a_j \,a_i \,, \qquad
 b_i \,b_j =b_j \,b_i \,, \\
{}& \quad  a_i \, a_i^{-1} = a_i^{-1} \, a_i =1 \,, \qquad
 b_i \, b_i^{-1} = b_i^{-1} \, b_i =1 \,,
\end{aligned}
\end{equation}
where $\delta_{ij}$ stands for the Kronecker symbol.

Let us remark that $Q_q(n)$ admits a representation
$\rho$ such that $\rho(a_i) = \exp (\alpha A_i)$
and $\rho(b_i) = \exp (\beta A^+_i)$, $q=e^{\alpha\beta}$,
where $A^+_i$ and $A_i$ are the creation and annihilation 
operators on the Fock space for $n$ non--interacting
particles, i.e. $[A_i,A^+_j]= 1\,\delta_{ij}$ and
\hbox{$[A^+_i,A^+_j]=[A_i,A_j]=0$}.
Replacing the quantum exponential function (\ref{Sx}) 
with the modular quantum dilogarithm \cite{Fa1} which 
has no divergences at $|q|=1$, one can also consider 
a representation $\rho'$ of $Q_q(n)$ where the generators 
are realised by unitary operators on the Hilbert space
$L^2({\mathbb R}^n)$:  
$\rho'(a_i) = \exp(\alpha \partial_{x_i})$, 
$\rho'(b_i) = \exp(\sqrt{-1}\, \beta x_i)$,
$q=e^{\sqrt{-1}\, \alpha\beta}$, $\alpha, \beta \in \mathbb R$.

Let $\phi : D_{q,r}(2) \to Q_{q}(1)$ be 
the homomorphism such that 
$\phi  (x_{11}) = a$,
$\phi  (x_{12}) = b$,
$\phi  (x_{22}) = a^{-r}$.
It can be promoted to two  
homomorphisms  
$\phi^{(k)} : D_{q,r}(3) \to Q_{q}(1)$
that correspond to two embeddings 
$SL^+_r(2) \subset SL^+_r(3)$ of the classical groups, namely:
\begin{align}
\label{phi1} {}&
 \phi^{(1)} (x_{11}) = a \,, \quad 
 \phi^{(1)} (x_{12}) = b \,, \quad
 \phi^{(1)} (x_{22}) = a^{-r} \,, \quad
 \phi^{(1)} (x_{k3}) = \delta_{k3} \,, \\[1mm]
\label{phi2} {}& 
\phi^{(2)} (x_{22}) = a \,, \quad 
 \phi^{(2)} (x_{23}) = b \,, \quad
 \phi^{(2)} (x_{33}) = a^{-r} \,, \quad
 \phi^{(2)} (x_{1k}) = \delta_{k1} \,.
\end{align}

Given $\lambda \in {\mathbb C}$, choose $\eta$ 
such that $\eta^{r^2+r+1} = \lambda^{-r}$. 
Consider two matrices:
\begin{align}
\label{B121}
{} 
T &= \eta \left( \begin{matrix}
 a_1 & b_1 & \\ & a^{-r}_1 & \\ && 1
\end{matrix} \right)
\left( \begin{matrix}
 1 && \\ &  a_2 & b_2 \\ && a^{-r}_2 
\end{matrix} \right)
\left( \begin{matrix}
 1 &  & \\ & \lambda & \\ && 1
\end{matrix} \right)
\left( \begin{matrix}
 a_3 & b_3 & \\ & a^{-r}_3 & \\ && 1
\end{matrix} \right) , \\[1mm]
\label{B212}
{}
 \tilde{T} &=
 \eta \left( \begin{matrix}
 1 && \\ &  a_3 & b_3 \\ && a^{-r}_3 
\end{matrix} \right)
\left( \begin{matrix}
 1 &  & \\ & \lambda & \\ && 1
\end{matrix} \right)
 \left( \begin{matrix}
 a_2 & b_2 & \\ & a^{-r}_2 & \\ && 1
\end{matrix} \right)
\left( \begin{matrix}
 1 && \\ &  a_1 & b_1 \\ && a^{-r}_1 
\end{matrix} \right)  .
\end{align}

\begin{lem}\label{BBt}
$T$ and $\tilde{T}$ satisfy relation (\ref{RGG})
with R--matrix $R_{q,q^r}^{(3)}$.
For $T$ and $\tilde{T}$,
the value of $C_0$ given by (\ref{c0}) 
is equal to one.
\end{lem}

Thus, formulae (\ref{B121}) and (\ref{B212}) define
(but not uniquely)
homomorphisms 
$\phi^{(121)}, \phi^{(212)}:  D_{q,r}(3)
 \to Q_{q}(3)$
such that $T=\phi^{(121)} \bigl(X(3)\bigr)$,
$\tilde{T}=\phi^{(212)} \bigl(X(3)\bigr)$.
Under these homomorphisms,
the images of the central element $C_1$ are 
monomials in $a_i$ and~$b_i$. For us, it will
be convenient to use 
\begin{align}
\label{cct1}  
{}& w \equiv q \eta^{-1}
  \phi^{(121)} \bigl( C_1^{-1} \bigr)
 = a_1^{r} a_2^{-r}  b_1 b_3^{-1} a_3 \,, \\
\label{cct2} 
{}& \tilde{w} \equiv q \eta^{-1}
  \phi^{(212)} \bigl( C_1^{-1} \bigr)
 =  a_3^{-r} b_1 b_3^{-1} a_1^{-1} a_2   \,.
\end{align}  
Observe that
\begin{equation}\label{wwff} 
     \tilde{w} = 
    (a_1^{-1} a_2 a_3^{-1})^{r+1} \, w 
  = \eta \lambda^{-1}
  \phi^{(121)}\bigl(x_{11}^{-1} x_{22} x_{33}^{-1} \bigr) \, w 
= q \lambda^{-1} 
\phi^{(121)}\bigl(x_{11}^{-1} x_{22} x_{33}^{-1} C_1^{-1 }\bigr)
  \,.
\end{equation}
Whence it follows that 
\begin{equation}\label{wwt} 
   w  \, \tilde{w} = \tilde{w} \, w \,.
\end{equation}

Now our aim is to construct
solutions to the equation  
\begin{equation}\label{RBBB} 
 \SR(\lambda) \  T    =
 \tilde{T}   \ \SR(\lambda) \,. 
\end{equation}
First, we consider an automorphism $\CR$ of 
$Q_q(3)$ such that
\begin{equation}\label{RBBBa} 
 \CR( T )   =
 \tilde{T}   \,, 
\end{equation}
where the action of $\CR$ on the l.h.s. is entry--wise.

\begin{prop}\label{Raut}
(i) Equation (\ref{RBBBa}) is equivalent to the following
set of relations:
\begin{align}
\label{aa'}  
{}& \CR(a_1 a_3) = a_2 \,, \qquad
\CR(a_2)  = a_1 a_3 \,, \qquad
\CR(b_1 b_2) = b_1 b_2 \,, \\
\label{ab'}
{}& \CR(a_1^{-r} b_2) = 
  a_1^{-r} b_3\, (1 + \lambda\, w)\,, \qquad
  \CR(a_1  b_3) = (1 + \lambda\, w)^{-1} a_1 b_2 \,.
\end{align}  
(ii)  $\CR$ acts on  $w$ and $\tilde{w}$ as follows:
\begin{equation}\label{Rwwg} 
   \CR(w) =\tilde{w} \,, \qquad
  \CR(\tilde{w}) = w \,.
\end{equation}
\end{prop}

Setting $\lambda=0$ in (\ref{aa'})--(\ref{ab'}),
we obtain another homomorphism.
Clearly, since we have five relations (\ref{aa'})--(\ref{ab'})
for six generators, this homomorphism is not defined
uniquely.

\begin{prop}\label{FF}
(i) If $\CF$ is a homomorphism satisfying (\ref{aa'})--(\ref{ab'})
for $\lambda=0$ 
and $\CF$ is an involution for all $r \neq -1$, then it is 
determined uniquely and is given by
\begin{align}
\label{Fa} &
 \CF ( a_1 )   = a_1  \,, \qquad
 \CF( a_2 )   = a_1   a_3 \,, \qquad
 \CF( a_3 )   = a_1^{-1} a_2 \,, \\ 
\label{Fb}  &
 \CF( b_1 )    = b_1 b_2 b_3^{-1} \,, \qquad
 \CF( b_2 )  = b_3 \,, \qquad
 \CF( b_3 ) = b_2 \,.
\end{align} 
(ii) $\CF$ acts on  $w$ and $\tilde{w}$ as follows:
\begin{equation}\label{Fwwg} 
  \CF(w) =\tilde{w} \,, \qquad
  \CF(\tilde{w}) = w \,.
\end{equation}
\end{prop}

Let $\CF_{ijk}$, $1 \leq i < j<k \leq 6$, 
denote the automorphism of $Q_q(6)$
which acts nontrivially only on the $a$'s and $b$'s
with labels $i,j,k$ 
and its action on these variables is given by
(\ref{Fa})--(\ref{Fb}), where 1 is substituted by $i$,
2 by $j$, and 3 by~$k$. It is known 
(see Proposition~2.1 in \cite{KV2} and Lemma~2.13
in~\cite{BV}) that $\CF$ satisfies the tetrahedron equation,
\begin{equation}\label{Ftetr}
 \CF_{123} \circ \CF_{145} \circ \CF_{246} \circ \CF_{356} =
 \CF_{356} \circ \CF_{246} \circ \CF_{145} \circ \CF_{123} \,.
\end{equation}
This statement can be verified by computing 
the action of both sides of (\ref{Ftetr}) on the monomial 
$M=a_1^{\alpha_1} b_1^{\beta_1} \ldots a_6^{\alpha_6} b_6^{\beta_6}$.

Thus, since $\CF$ is an involution and a solution
to the tetrahedron equation, 
it is a counterpart of the permutation ${\mathcal P}^{(n)}$
in the Yang--Baxter case. Therefore, 
by analogy with the Yang--Baxter case, one can
look for solutions to (\ref{aa'})--(\ref{ab'})
of the form $\CR = \check{\CR} \circ \CF$.
Using (\ref{aa'})--(\ref{ab'}) and 
(\ref{Fa})--(\ref{Fb}), we infer that $\check{\CR}$
is an automorphism of $Q_q(3)$ such that 
\begin{align}
\label{aaa'}  
{}& \check{\CR}(a_1 a_3) = a_1 a_3 \,, \qquad
\check{\CR}(a_2)  = a_2 \,, \qquad
\check{\CR}(b_1 b_2) = b_1 b_2 \,, \\
\label{aab'}
{}& \check{\CR}(a_1^{-r} b_3) = 
  a_1^{-r} b_3\, (1 + \lambda\, w)\,, \qquad
  \check{\CR}(a_1  b_2) = (1 + \lambda\, w)^{-1} a_1 b_2 \,.
\end{align}   
Note that  $w$ and $\tilde{w}$ 
are invariant under the action of~$\check{\CR}$, 
\begin{equation}\label{CRwwg} 
   \qquad \check{\CR}(w) = w \,, \qquad
  \check{\CR}(\tilde{w}) = \tilde{w} \,.
\end{equation}

Let us look for an invertible element 
$\check{\SR}(\lambda) \in Q_{q}(3)$
such that 
\begin{equation}\label{Req} 
\check{\SR}(\lambda) \, x = 
\check{\CR}(x) \, \check{\SR}(\lambda)
\end{equation}
for every $x \in Q_{q}(3)$.
{}From (\ref{aaa'}) and (\ref{CRwwg}) it follows that
$\check{\SR}(\lambda)$ commutes with $a_1 a_3$, $a_2$, $b_1b_2$, $w$, 
and $\tilde{w}$.
It is not hard to check that the subalgebra in $Q_{q}(3)$
of elements commuting
with this set is generated by $w$ and $\tilde{w}$. Therefore, 
$\check{\SR}(\lambda)$ is a function 
of $w$ and $\tilde{w}$ only.
 
\begin{prop}\label{RWW}
Every invertible 
$\check{\SR}(\lambda) \in {\mathbb C}[[w,\tilde{w}]]$ 
which satisfies (\ref{Req}) is given by   
\begin{equation}\label{CRsol}
 \check{\SR}(\lambda) = f_{q,\lambda}(\tilde{w}) \,
 	\dil{\lambda\,w}{q^{r+1}} \,, 
\end{equation}
where $\dil{t}{q}$ was defined in (\ref{Sx})
and $f_{q,\lambda}(t)$ is an arbitrary invertible 
series from ${\mathbb C}[[t]]$.
\end{prop} 

Note that all $a_i$ and
$b_i$ $q$--commute with $w$ and $\tilde{w}$.
Therefore, if $M$ is a monomial in $a_i$ and $b_i$,
we have $\check{\CR}(M) = M\, h_1(w) \, h_2(\tilde{w})$,
where $h_1(t), h_2(t) \in {\mathbb C}[[t]]$.

Proposition~\ref{RWW} implies that we have constructed
solutions to equation (\ref{RBBB}) of the
following form:
\begin{equation}\label{RFs}
 \SR(\lambda) = 
 f_{q,\lambda}(\tilde{w})\,
 	\dil{\lambda\,w}{q^{r+1}} \,  \SF
 = \SF \,
 f_{q,\lambda}(w)\,
 	\dil{\lambda\,\tilde{w}}{q^{r+1}}   \,,
\end{equation}
where the formal element $\SF$ is an involution and 
$\SF\,x=\CF(x) \, \SF$ for every
$x \in Q_q(3)$. Note that the later relation cannot be 
resolved if $\SF$ is assumed to be a formal power
(or a formal Laurent) series in elements from~$Q_q(3)$.
Therefore, we will treat $\SF$ as the generator of
an outer automorphism.

Below we will need some generalization of
(\ref{RFs}) constructed by the same method.
For this purpose,  
we consider a family of automorphisms $\psi_s$,
$s \in {\mathbb Z}$, of $Q_q(3)$
defined on the generators as follows:
\begin{align}
\label{psia}  
{}&  \psi_s(a_1)=a_1 \,, \qquad 
 \psi_s(a_2)=a_2 \,, \qquad 
 \psi_s(a_3)=a_3    \,,\\
\label{psib}  
 \psi_s(b_1) & =   a_2^{s} a_3^{s} \, b_1 \,,  \qquad 
 \psi_s(b_2)= a_1^{s} \, b_2 \,, \qquad 
 \psi_s(b_3)= a_1^{s}  \,b_3    \,.
\end{align}   
Clearly, we have $\psi_s^{-1}=\psi_{-s}$.
It is straightforward to verify that $\psi_s$
commutes with $\CF$ given by~(\ref{Fa})--(\ref{Fb}),
\begin{align}\label{psiF}  
 \psi_s \circ \CF = \CF \circ \psi_s \,.
\end{align}  

Consider $T_s \equiv \psi_s(T)$ and 
$\tilde{T}_s \equiv \psi_s(\tilde{T})$,
that is matrices given by (\ref{B121}) and
(\ref{B212}), where $b_i$ are replaced with 
$\psi_s(b_i)$. Since $\psi_s$ is an automorphism of 
$Q_q(3)$,
Lemma~\ref{BBt} applies to 
$T_s$ and $\tilde{T}_s$ as well.
Therefore, we can repeat the construction 
of this section in order to find solutions
to the equation 
\begin{equation}\label{RBBBs} 
 \SR^{(s)}(\lambda) \ T_s    =
 \tilde{T}_s   \ \SR^{(s)}(\lambda) 
\end{equation} 
of the form $\SR^{(s)}(\lambda) = 
\check{\SR}^{(s)} (\lambda) \, \SF$,
where $\SF$ is the same as above. 
It is clear that 
$\CR_s \equiv \psi_s \circ \CR \circ \psi_{-s}$
is an automorphism of 
$Q_q(3)$
such that 
$\CR_s( T_s ) = \tilde{T}_s$.
Taking (\ref{psiF}) into account, we deduce
that $\CR_s= \check{\CR}_s \circ \CF$, where
$\check{\CR}_s = \psi_s \circ \check{\CR} \circ \psi_{-s}$.
{}From the latter relation along with  (\ref{Req}) and (\ref{CRsol})
we infer that $\check{\CR}_s$ corresponds in the
sense of Proposition~\ref{RWW} to 
$\psi_s\bigl(\check{\SR}(\lambda)\bigr)$.
Therefore, we draw the following conclusion:

\begin{thm}\label{RFG}
Equation (\ref{RBBBs})  has solutions of 
the following form:
\begin{equation}\label{RFsss}
 \SR^{(s)}(\lambda) = 
 f_{q,\lambda}\bigl(\tilde{w}^{(s)}\bigr)\,
 	\dil{\lambda\,w^{(s)}}{q^{r+1}} \,  \SF
 = \SF \ 
 f_{q,\lambda}\bigl(w^{(s)}\bigr)\,
 	\dil{\lambda\,\tilde{w}^{(s)}}{q^{r+1}}   \,,
\end{equation}
where $f_{q,\lambda}(t)$ is an arbitrary invertible 
series from ${\mathbb C}[[t]]$ and
\begin{align}\label{psiww}  
 w^{(s)} \equiv
 \psi_s (w) = a_1^{r-s} a_2^{s-r}  b_1 b_3^{-1} a_3^{s+1} 
  \,, \qquad
  \tilde{w}^{(s)} \equiv 
  \psi_s (\tilde{w}) = 
   a_3^{s-r} b_1 b_3^{-1} a_1^{-s-1} a_2^{s+1} \,.
\end{align}   
\end{thm}
Note that 
setting $s=0$ in (\ref{RFsss}) we recover (\ref{RFs}).
 
Applying the homomorphism $\psi_s$ to (\ref{wwt}),
we obtain 
\begin{equation}\label{wwt2} 
   w^{(s)}  \, \tilde{w}^{(s)} = \tilde{w}^{(s)} \, w^{(s)} \,.
\end{equation}
In addition, it follows from (\ref{Fwwg}) and (\ref{psiF})
that $\SF \, w^{(s)} = \tilde{w}^{(s)} \, \SF$.

\section{The tetrahedron equation}
In order to discuss solutions to the tetrahedron equation we 
introduce in the standard way tensor copies of
$\SR^{(s)}(\lambda)$, namely
\begin{equation}\label{RFFs}
 \SR^{(s)}_{ijk}(\lambda) = 
 f_{q,\lambda}(\tilde{w}^{(s)}_{ijk})\,
 	\dil{\lambda\,w^{(s)}_{ijk}}{q^{r+1}} \,  \SF_{ijk}
 = \SF_{ijk} \,
 f_{q,\lambda}(w^{(s)}_{ijk})\,
 	\dil{\lambda\,\tilde{w}^{(s)}_{ijk}}{q^{r+1}}   \,,
\end{equation}
where $1 \,{\leq}\, i \,{<}\, j \,{<}\, k \,{\leq}\, 6$. 
Explicitly, the arguments in (\ref{RFFs}) are 
the following elements of $Q_q(6)$:
\begin{align}\label{www}  
 w^{(s)}_{ijk}  =
  a_i^{r-s} a_j^{s-r}  b_i b_k^{-1} a_k^{s+1} 
  \,, \qquad
  \tilde{w}^{(s)}_{ijk} =
   a_k^{s-r} b_i b_k^{-1} a_i^{-s-1} a_j^{s+1} \,.
\end{align}   

The key feature of equation (\ref{RBBBs}) is that  
 $\SR^{(s)}(\lambda)$ is closely related  to 
 the transformation of the word 
$s_1 s_2 s_1$ into $s_2 s_1 s_2$. 
To make it explicit, let us rewrite (\ref{RBBBs}) as
follows:
\begin{equation}\label{RBD} 
  \SR^{(s)}_{ijk}(\lambda) \
  B^{(1)}_i B^{(2)}_j  D^{(12)}_\lambda B^{(1)}_k  =
 B^{(2)}_k D^{(12)}_\lambda  B^{(1)}_j   B^{(2)}_i \ 
 \SR^{(s)}_{ijk}(\lambda) \,,
\end{equation}
where the notations are self--evident (we omit index $s$
in the matrix factors).

For the Weyl group of $A_3$, the two reduced forms
of the longest element that are most distant from
each other are $s_1 s_2 s_1 s_3 s_2 s_1$
and $s_3 s_2 s_1 s_3 s_2 s_3$.
Therefore,  
we consider two $4\,{\times}\,4$ matrices
with entries in $Q_q(6)$
corresponding to these words (we again omit index $s$
in the matrix factors):
\begin{align} 
{}&
 A_s=B^{(1)}_1 B^{(2)}_2  D^{(12)}_\lambda B^{(1)}_3 
 B^{(3)}_4 D^{(23)}_{\lambda\nu} B^{(2)}_5 D^{(12)}_\mu  B^{(1)}_6 
  \,, \\
{}&
 \tilde{A}_s = B^{(3)}_6 D^{(23)}_{\nu} B^{(2)}_5  
 D^{(12)}_{\lambda\mu} B^{(1)}_4 
 B^{(3)}_3  D^{(23)}_{\lambda} B^{(2)}_2  B^{(3)}_1    \,. 
\end{align}
Note that $B^{(1)}_k$ and $B^{(3)}_j $ commute.

Recall that the word $s_1 s_2 s_1 s_3 s_2 s_1$ can be transformed
into $s_3 s_2 s_1 s_3 s_2 s_3$ along two different paths
composed of four local steps, see~(\ref{We123}). 
By (\ref{RBD}), every local braid transformation
$s_l s_{l+1} s_l = s_{l+1} s_l s_{l+1}$ in (\ref{We123})
corresponds to  $\SR^{(s)}_{ijk}(\lambda)$ 
for some $i,j,k$. Therefore, 
one can expect that 
$A_s$ is transformed into~$\tilde{A}_s$ 
by products of four~$\SR$'s.

\begin{prop}\label{AAt}
(i) $A_s$ and $\tilde{A}_s$ 
satisfy the RTT relation (\ref{RGG})
with R--matrix $R_{q,q^r}^{(4)}$. \\
(ii) The following relations hold:
\begin{align}
\label{RBB4a} 
 \SM \  A_s= 
 \tilde{A}_s \  
 \SM \,,  \qquad\quad 
 \SM'  \  A_s = 
 \tilde{A}_s \ \SM'   ,
\end{align} 
where 
\begin{align}
\label{MM} 
\SM = \SR^{(s)}_{123} (\lambda) \SR^{(s)}_{145} (\lambda\mu) 
 \SR^{(s)}_{246} (\nu) \SR^{(s)}_{356}(\mu)  , \qquad 
 \SM' = \SR^{(s)}_{356}(\nu) \SR^{(s)}_{246}(\mu) 
 \SR^{(s)}_{145}(\lambda\nu)  
 \SR^{(s)}_{123}(\lambda) \,. 
\end{align} 
\end{prop}

Part (i) of this Proposition  implies that
$A_s$ and $\tilde{A}_s$ (multiplied by suitable scalars) 
define two homomorphisms 
{}from $D_{q,r}(4)$ to 
$Q_q(6)$.
However the center of $D_{q,r}(4)$ is not scalar.
Therefore, relations (\ref{RBB4a}) do not imply
that $\SM = \SM'$, which would mean that 
$\SR^{(s)}(\lambda)$ satisfies the
tetrahedron equation
\begin{equation}\label{Rtetr2}
 R_{123} (\lambda) \, R_{145} (\lambda\mu) \, 
  R_{246} (\nu) \, R_{356}(\mu) =
  R_{356}(\nu) \, R_{246}(\mu) \, 
  R_{145}(\lambda\nu) \, R_{123}(\lambda) \,,
\end{equation}
but only that $\SM^{-1} \SM'$ 
commutes with all the entries of~$A_s$.

In order to find out whether $\SR^{(s)}(\lambda)$
given by (\ref{RFsss}) satisfies the
tetrahedron equation, we consider a more general ansatz.

\begin{thm}\label{RS} 
(i) Let $a_i$, $b_i$, $i=1,{\ldots},6$ satisfy 
relations (\ref{B222}) and 
$g_{q,\lambda}(t) \in {\mathbb C}[[t]]$ be a formal power
series such that $g_{q,\lambda}(0)=1$. Then
 \begin{equation}\label{Rgl}
 R_{123}(\lambda) = \SF_{123} \ 
 g_{q,\lambda}( q^{-\alpha_1 \beta_1}
 a_1^{\alpha_1} a_2^{\alpha_2} a_3^{\alpha_3} 
  b_1^{\beta_1} b_2^{\beta_2} b_3^{\beta_3})  ,
\end{equation}
where $\alpha_i, \beta_i \in {\mathbb Z}$,
satisfies the tetrahedron equation (\ref{Rtetr2})
provided that  
\begin{equation}\label{xy}
 \alpha_1 + \alpha_2 = 0 \,, \qquad \beta_2 =0 \,, \qquad
 \beta_1 + \beta_3 = 0 \,,
\end{equation}
and 
\begin{equation}\label{ru}
 g_{q,\lambda}(t)= 
 \dil{Q^{\frac{m(m-1)}{2}} \lambda^k \, t^m}{Q^{m^2}}  \,,
\end{equation}
where $m\in {\mathbb Z}_+$, $k\in {\mathbb Z}$, and 
\begin{equation}\label{Qru}
 Q = q^{-(\alpha_1 + \alpha_3)\beta_1} \,.
\end{equation}
(ii) Let 
$g_{q,\lambda}(t), f_{q,\lambda}(t) \in {\mathbb C}[[t]]$ 
be formal power series such that 
$g_{q,\lambda}(0)=f_{q,\lambda}(0)=1$. 
Then
 \begin{equation}\label{Rgfl}
 R_{123}(\lambda) = \SF_{123} \ 
 g_{q,\lambda}( q^{-\alpha \beta}
 a_1^{\alpha} a_2^{-\alpha} a_3^{\delta} \,
  b_1^{\beta} b_3^{-\beta})  \
   f_{q,\lambda}(q^{-\hat{\alpha} \hat{\beta} }
  a_1^{\hat{\alpha}} a_2^{-\hat{\alpha}} 
  a_3^{\hat{\delta}} \,
  b_1^{\hat{\beta}} b_3^{-\hat{\beta}}) ,
\end{equation}
where $\alpha, \beta, \delta, \hat{\alpha},
\hat{\beta}, \hat{\delta} \in {\mathbb Z}$,
satisfies the tetrahedron equation (\ref{Rtetr2})
provided that  
\begin{equation}\label{xy22}
 \alpha \, \hat{\beta} + \hat{\delta} \, \beta =0 \,, \qquad
  \hat{\alpha} \, \beta + \delta \, \hat{\beta}   =0 \,,
\end{equation}
and 
\begin{equation}\label{ruru}
 g_{q,\lambda}(t)= 
 \dil{Q^{\frac{m(m-1)}{2}} \lambda^k \, t^m}{Q^{m^2}}  \,,
 \qquad
 f_{q,\lambda}(t)=   
 \dil{ \widehat{Q}^{\frac{\widehat{m}(\widehat{m}-1)}{2}} 
  \lambda^{l}  \, t^{\widehat{m}}}{ {\widehat{Q}}^{\widehat{m}^2}} ,
\end{equation}
where $m,\widehat{m}\in {\mathbb Z}_+$, $k,l \in {\mathbb Z}$, 
and 
\begin{equation}\label{QQru}
 Q = q^{-(\alpha  + \delta)\beta} \,, \qquad
 \hat{Q} = q^{-(\hat{\alpha}  + \hat{\delta}) \hat{\beta}} \,.
\end{equation}
\end{thm}

Let us remark that solutions to the tetrahedron equations 
of the type (\ref{Rgl}) were found earlier by various 
authors (usually for $\beta_1=m=1$). The case $\alpha_1=0$, 
 $\alpha_3 =-1$ 
was considered in \cite{Se1}, the 
case $\alpha_1=\alpha_3 =-1$ 
in \cite{KV2}, the case 
$\alpha_3 =-1-\alpha_1$   in~\cite{BV}.

It is a non--trivial fact that a solution of the type
(\ref{Rgl}) admits the second factor, $f_{q,\lambda}$.
Note that  conditions (\ref{xy22}) imply that
the arguments of $g_{q,\lambda}$ and $f_{q,\lambda}$
in (\ref{Rgfl}) commute. Furthermore,
(\ref{xy22}) and (\ref{QQru}) imply that the corresponding
$q$--parameters are related as follows:
\begin{equation}\label{QQQ}
 Q^{\widehat{\beta}^2 } = {\widehat{Q}}^{- \beta^2} \,.
\end{equation}

Observe that the arguments of $g_{q,\lambda}$ and 
$f_{q,\lambda}$ in (\ref{Rgfl})
coincide with $\tilde{w}^{(s)}$ and $q^{-(r+1)}w^{(s)}$,
respectively, if we set
\begin{equation}\label{wwab}
 \alpha = - \hat{\delta} = -s-1 \,, \qquad 
 \hat{\alpha} = - \delta = r -s \,, \qquad
 \beta = \hat{\beta} = 1 \,.
\end{equation}
By (\ref{QQru}), this corresponds to
$Q=q^{r+1}$ and $\widehat{Q}=q^{-(r+1)}$.
Therefore, Theorem~\ref{RS} has the following 
corollaries concerning solutions to the 
intertwining relation~(\ref{RBBBs}):

\begin{cor} 
$\SR^{(s)}(\lambda)$ given by (\ref{RFFs})
satisfies the tetrahedron equation (\ref{Rtetr2})
provided that invertible formal power series 
$f_{q,\lambda}(t) \in {\mathbb C}[[t]]$ 
either is constant or it is given by
\begin{equation}\label{ru3}
f_{q,\lambda}(t)= f_{q,\lambda}(0) \,
\dil{ q^{-(r+1)\frac{\widehat{m}(\widehat{m}+1)}{2}} 
 \lambda^{l}  \, t^{\widehat{m}}}{ q^{-(r+1)\widehat{m}^2}} \,,  
\end{equation}
where $\widehat{m}\in {\mathbb Z}_+$ and $l \in {\mathbb Z}$.
\end{cor}

\begin{cor}
If $\SR^{(s)}(\lambda)$ given by (\ref{RFFs})
is an involution, i.e.,
$\bigl(\SR^{(s)}(\lambda)\bigr)^2 = 1$,
then it satisfies the tetrahedron equation~(\ref{Rtetr2}).
\end{cor}

Indeed, multiplying both expressions for
$\SR^{(s)}(\lambda)$ in (\ref{RFFs}) and 
taking into account that $\SF^2=1$, we deduce
that $(\SR^{(s)}(\lambda))^2 = 1$ holds provided that
\begin{equation}\label{rrss}
 f_{q,\lambda} \bigl(w^{(s)}\bigr) \, 
 \dil{\lambda\,w^{(s)}}{q^{r+1}} =
 \Bigl( f_{q,\lambda} \bigl(\tilde{w}^{(s)} \bigr)  
 \dil{\lambda\,\tilde{w}^{(s)}}{q^{r+1}} \Bigr)^{-1} \,.
\end{equation}
Since $w^{(s)}$ and $\tilde{w}^{(s)}$ are algebraically independent, 
the formal power series on the l.h.s. and r.h.s. of
(\ref{rrss}) can coincide only if 
$f_{q,\lambda}(t)= f_{q,\lambda}(0)/\dil{\lambda\,t}{q^{r+1}}$,
where $f_{q,\lambda}(0)=\pm 1$. 
In this case, $\SR^{(s)}(\lambda)$ is similar to~$\SF$ 
\begin{equation}\label{RFss}
 \SR^{(s)}(\lambda) = 
 \pm \dil{\lambda\,w^{(s)}}{q^{r+1}}
   \SF \,  
  \frac{1}{\vphantom{\Bigm|} \dil{\lambda\,w^{(s)}}{q^{r+1}}}
  = 
 \pm 
 \frac{1}{\vphantom{\Bigm|} \dil{\lambda\,\tilde{w}^{(s)}}{q^{r+1}}} 
  \ \SF \,  \dil{\lambda\,\tilde{w}^{(s)}}{q^{r+1}} .
\end{equation}
Note that, by (\ref{Sxi}), we have 
$f_{q,\lambda}(t) =  
f_{q,\lambda}(0)/\dil{\lambda\,t}{q^{r+1}} =
f_{q,\lambda}(0)\, \dil{q^{-r-1}\lambda\,t}{q^{-(r+1)}}$
which coincides with (\ref{ru3}) for $\hat{m}=l=1$.
Therefore, by Corollary~1, such $\SR^{(s)}(\lambda)$
satisfies the tetrahedron equation.

\appendix
\section{Appendix}

\begin{proof}[\bf Proof of Lemma~\ref{pentdil}] 
For the series 
$S_{q,\lambda}(t)=S_0(\lambda) +S_1(\lambda) t+ \ldots$
(we omit the dependence on $q$), we have 
$S_0(\lambda) \neq 0$ (since otherwise the monomials of 
the least powers on the l.h.s. and on the r.h.s. of (\ref{pent}) 
do not match). Moreover, we have
$S_0(\lambda)S_0(\lambda\mu)S_0(\mu)=S_0(\mu)S_0(\lambda)$,
whence we conclude that $S_0(\lambda) =1$.
 Since relations (\ref{XYZ}) 
are homogeneous in $\SX$ and $\SY$, the coefficients at $\SY^k$
on the l.h.s. and on the r.h.s. of (\ref{pent}) 
must match for every~$k$. Let $m$ be the 
smallest positive integer such that $S_m(\lambda) \neq 0$. 
Matching 
the coefficients at $\SY^m$ on the both sides of (\ref{pent}), 
we infer that  
\begin{equation}\label{Seq}
  S_{q,\lambda}(\SX) \, 
  \bigl(S_m(\mu)\, \SY^m +  S_m(\lambda\mu) \, 
  (\SX \, \SY)^m \bigr) = 
  S_m(\mu)\, \SY^m \, S_{q,\lambda}(\SX) \,.
\end{equation}
Taking into account that $\SY \SX = q \, \SX \SY$,
we can rewrite (\ref{Seq}) as follows:
\begin{equation}\label{Seq2}
   S_{q,\lambda}(\SX) \, 
  \bigl(1 +
  q^{\frac{m(m-1)}{2}} \frac{S_m(\lambda\mu)}{S_m(\mu)}\, \SX^m \bigr)
  = S_{q,\lambda}(q^m \SX)  \,.
\end{equation}
It follows that $\frac{S_m(\lambda\mu)}{S_m(\mu)}$
does not depend on $\mu$ and hence
$S_m(\lambda) \sim \lambda^k$, $k\in \mathbb Z$.
Thus, (\ref{Seq2}) acquires the form
\begin{equation}\label{Seq3}
  S_{q,\lambda}(q^m \SX) = S_{q,\lambda}(\SX) \, 
  \bigl(1 + q^{\frac{m(m-1)}{2}} \lambda^k \, \SX^m \bigr) \,.
\end{equation}
After a change of variables,
$\hat{\SX}=q^{\frac{m(m-1)}{2} } \lambda^k \, \SX^m$,
$\hat{q} = q^{m^2}$,
functional equation (\ref{Seq3}) turns into 
equation (\ref{fundil}). Whence we 
deduce that 
$S_{q,\lambda}(t)=\dil{q^{m(m-1)/2} \lambda^k \, t^m}{q^{m^2}}$. 

By the same reasoning as above, we deduce that if the
series
$U_{q,\lambda}(t) = U_0(\lambda) + U_1(\lambda)t + \ldots$
is not vanishing, then $U_0(\lambda)$ does not vanish
and does not depend on~$\lambda$. Therefore, setting
$\SZ=0$ in (\ref{tetr}), we conclude  that
$U'_{q,\lambda}(t) \equiv U_{q,\lambda}(t)/U_0$ satisfies (\ref{pent}) 
and hence 
$U_{q,\lambda}(t)=U_0 \,\dil{q^{m(m-1)/2} \lambda^k \, t^m}{q^{m^2}}$
for some $m$ and~$k$. In order to verify that 
relation (\ref{tetr}) holds for such a series one
has to make the same change of variables as above
and compare the resulting relation with identity~(\ref{qexp4}).
\end{proof}

\begin{proof}[\bf Proof of Lemma~\ref{Cdet}]
Consider the element
\begin{align}\label{c00}
 C_0 = \prod_{i=1}^n (x_{ii})^{\alpha_i}  \,,
\end{align}  
where $\alpha_i$ are unknown integers. The element $x_{jm}$, $j<m$,
commutes  with $x_{ii}$ non--trivially only if $i \in [j,m]$.
Specifically, we have $x_{jm} x_{ii} = \theta\, x_{ii} x_{jm}$ ,
where $\theta = q^{-1}$ if $i=j$, $\theta = q^{r}$ if $i=m$,
and $\theta = q^{r-1}$ if $j<i<m$. Using these relations,
we deduce that $x_{jm}$ commutes with $C_0$ provided that
$\alpha$'s satisfy the following relations
for all pairs $j,m$, where $j<m$:
\begin{align}\label{ac0}
 -\alpha_j + r\,\alpha_m + (r-1) \sum_{j<i<m} \alpha_i = 0   \,.
\end{align}
It is easy to verify that (\ref{ac0}) holds 
if we choose  $\alpha_i = r^{n-i}$.
\end{proof} 

\begin{proof}[\bf Proof of Lemma~\ref{BBt}]
By the standard argument, if two
matrices $X$ and $X'$ satisfy (\ref{RGG}) and
entries of $X$ commute with these of $X'$,
then (\ref{RGG})  holds for $XX'$ as well.
Every matrix factor in $T$ and $\tilde{T}$
is of the form $\phi\bigl(X(3)\bigr)$,
where $\phi$ is either one of the homomorphisms
(\ref{phi1})--(\ref{phi2}) or a trivial 
homomorphism sending $x_{ii}$ to $\lambda^{\delta_{i2}}$
and vanishing on $x_{ij}$ if $i\neq j$.
Therefore, every matrix factor
in $T$ and $\tilde{T}$ satisfies relation
(\ref{RGG}) with $\check{R}_{q,q^r}^{(3)}$.
Taking into account that entries of 
different factors in $T$ commute, we conclude that
$T$ satisfies
(\ref{RGG}) with $\check{R}_{q,q^r}^{(3)}$.
The same holds for $\tilde{T}$.
The corresponding values of~$C_0$
are straightforward to compute.
\end{proof} 

\begin{proof}[\bf Proof of Proposition~\ref{Raut}]

By a simple computation, we find
\begin{align} 
\label{B121e}
{}& 
 T = \eta \left( \begin{matrix}
 a_1 a_3 & a_1 b_3 + \lambda\, b_1 a_2 a^{-r}_3 &  b_1 b_2 \\ 
 & \lambda\, a^{-r}_1 a_2 a^{-r}_3 & \, a^{-r}_1  b_2 \\ 
 &&  a^{-r}_2
\end{matrix} \right) \,, \\ 
\label{B212e}
{}& \tilde{T} = \eta \left( \begin{matrix}
 a_2 &  a_1 b_2  &  b_1 b_2 \\ 
 & \lambda\, a_1 a^{-r}_2 a_3 & 
  b_3 a^{-r}_1 + \lambda\,  b_1 a^{-r}_2 a_3 \\ 
  &&  a^{-r}_1 a^{-r}_3
\end{matrix} \right) \,.
\end{align}
Now, comparing the diagonal entries and 
the upper right entries of $T$ and~$\tilde{T}$,
we obtain relations~(\ref{aa'}).
Further, comparing the entries $(2,3)$ and
taking formula (\ref{cct1}) into account, we obtain
the first relation in~(\ref{ab'}):
\begin{align}\label{rab} 
 \CR(a_1^{-r} b_2) =
 b_3 a^{-r}_1 + \lambda\,  b_1 a^{-r}_2 a_3
 =  a^{-r}_1 b_3 (1+ \lambda\, w)\,.
\end{align}
Comparing the entries $(1,2)$, we infer that
\begin{align}
 \CR(a_1 b_3) + \lambda\, \CR(b_1 a_2 a_3^{-r})
 =  a_1 b_2 \,.
\end{align}
Whence the second relation in~(\ref{ab'}) is derived
as follows:
\begin{align*} 
 \CR( a_1 b_3) & =  a_1 b_2  - \lambda\,
  \CR( b_1 b_2) \CR( (a^{-r}_1 b_2)^{-1})
  \CR( a^{-r}_1 a^{-r}_3 ) \CR(a_2)\\ 
 & \stackrel{(\ref{aa'}),(\ref{rab})}{=} 
  a_1 b_2  - \lambda\,
 b_1 b_2  (1+\lambda\, w)^{-1} 
 a_1^{r} a_2^{-r} b_3^{-1} a_3  a_1 \\
\nonumber   
 & \stackrel{(\ref{cct1})}{=}
   b_2 b_1 b_1^{-1} a_1 - \lambda\, 
 b_2 b_1 (1+\lambda\, w)^{-1} w b_1^{-1} a_1 \\  
& =   b_2 b_1 
 \bigl( 1 - \lambda\, (1+\lambda\, w)^{-1} w  \bigr)
 b_1^{-1} a_1 \\
 & =  b_2 b_1 (1+ \lambda\, w)^{-1} b_1^{-1} a_1  
   =  (1+ \lambda\, w)^{-1}  a_1 b_2.
\end{align*}
In the last line we used that $w$ commutes with
$b_1 b_2$.

Consider the element $z \equiv a_1^{1-r} b_2 b_3$.
Multiplying relations in (\ref{ab'}), we infer that
\begin{align}\label{crz} 
  \CR(z) = z \,.
\end{align} 
Note that
$w$ and $\tilde{w}$ can be written as follows:
\begin{align}\label{wwz} 
w=q a_2^{-r} z^{-1} (b_1 b_2) (a_1 a_3)  \,, \qquad
\tilde{w}= q (a_1 a_3)^{-r} z^{-1} (b_1 b_2) \, a_2\,.
\end{align}
Now, relations  (\ref{Rwwg})
are obvious if we take into account  
relations (\ref{aa'}) and~(\ref{crz}).
\end{proof}

\begin{proof}[\bf Proof of Proposition~\ref{FF}]
Setting $\lambda=0$ in (\ref{aa'})--(\ref{ab'}), we obtain
\begin{align}
\label{faa'}  
{}& \CF(a_1 a_3) = a_2 \,, \qquad
\CF(a_2)  = a_1 a_3 \,, \qquad
\CF(b_1 b_2) = b_1 b_2 \,,  \\
\label{fab'}
{}&\qquad  \CF(a_1^{-r} b_2) = 
  a_1^{-r} b_3\,, \qquad
  \CF(a_1  b_3) =  a_1 b_2 \,.
\end{align}

It follows from (\ref{fab'}) that 
$\CF(z)=z$. Therefore, formulae (\ref{wwz})
and (\ref{faa'}) imply that relations (\ref{Fwwg})
hold (even if $\CF$ is not an involution).

Consider $\CF(a_1)$. Relations (\ref{faa'})--(\ref{fab'})
imply that $\CF(a_1)$ commutes with 
$a_2$, $a_1 a_3$, $a_1^{-r} b_3$, $a_1 b_2$
and that $\CF(a_1) (b_1 b_2) =q (b_1 b_2) \CF(a_1)$.
It is easy to check that every monomial satisfying these 
relations is of the form 
$a_1^{1-k} a_2^{k} a_3^{-rk} b_1^k b_3^{-k}$,
$k \in {\mathbb Z}$. Comparing with (\ref{cct2}),
we infer that $\CF(a_1)=a_1\,h(\tilde{w})$,
where $h(t)$ is an unknown function.
By (\ref{Fwwg}), we have $\CF(\tilde{w})=w$.
Therefore, if
 $\CF$ is an involution, we have
$a_1 = \CF\bigl(\CF(a_1)\bigr)=
 \CF\bigl(a_1\,h(\tilde{w}) \bigr)=
a_1\,h(w)\,h(\tilde{w})$.
Whence, $h(w)\,h(\tilde{w})=1$ which implies that
$h(t)=\pm 1$ and $\CF(a_1)= \pm a_1$. 
For the minus sign, we have 
$\CF\bigl(\CF(a_1^{-r} b_2)\bigr)=
\CF(a_1^{-r} b_3) = \CF(a_1^{-r-1} a_1 b_3) =
(-1)^{-r-1} a_1^{-r} b_2$ and thus $\CF$ is
an involution only for odd~$r$.
Therefore, we have to choose the plus sign. 
Then, using $\CF(a_1)=a_1$,
 it is straightforward to derive
the rest of formulae in~(\ref{Fa})--(\ref{Fb}).
\end{proof}

\begin{proof}[\bf Proof of Proposition~\ref{RWW}] 

$\check{\SR}(\lambda)\equiv\check{\SR}_{\lambda}(w,\tilde{w})$ 
satisfies relations
(\ref{aaa'}) automatically because 
$a_1a_3$, $a_2$, and $b_1 b_2$
commute with both $w$ and $\tilde{w}$.
Further, we note that $a_1^{-r} b_3$ and $a_1 b_2$
commute with $\tilde{w}$ but have the following
non--trivial relations with $w$:
\begin{align}\label{wab}  
  w \, (a_1^{-r} b_3) = q^{r+1} \, (a_1^{-r} b_3) \, w \,,
  \qquad 
  w \, (a_1 b_2) = q^{-r-1} \, (a_1 b_2) \, w  \,.
\end{align}
Using these relations, it is not hard to verify that
both equations (\ref{aab'}) are equivalent  
to the following functional equation:
\begin{align}\label{rwq}  
   \check{\SR}_{\lambda}(q^{r+1} w , \tilde{w}) = 
   (1 + \lambda w) \, \check{\SR}_{\lambda}(w , \tilde{w})  \,.
\end{align}
This is the functional equation~(\ref{fundil}),
where $q$ is replaced with $q^{r+1}$ and 
$\tilde{w}$ plays the role of a constant.
Its solution in ${\mathbb C}[[w]]$ 
is unique up to a scalar multiple. The latter
can be a non--trivial formal power series in~$\tilde{w}$
whose coefficients can depend on $q$ and~$\lambda$. 
\end{proof}

\begin{proof}[\bf Proof of Proposition~\ref{AAt}]
Part (i). The same arguments apply as those 
that were used in the proof of Lemma~\ref{BBt}.\\
Part (ii). A verification is straightforward.
For instance, we have
\begin{align*}
{}&
 \SM' \ A_s = R^{(s)}_{356}(\nu) R^{(s)}_{246}(\mu) 
 R^{(s)}_{145}(\lambda\nu)  
 R^{(s)}_{123}(\lambda) \ A_s  \\
{}& =
 R^{(s)}_{356}(\nu) R^{(s)}_{246}(\mu) 
 R^{(s)}_{145}(\lambda\nu)  
 R^{(s)}_{123}(\lambda) \
 \underline{
 B^{(1)}_1 B^{(2)}_2  D^{(12)}_\lambda B^{(1)}_3 }
 B^{(3)}_4 D^{(23)}_{\lambda\nu} B^{(2)}_5 D^{(12)}_\mu  
 B^{(1)}_6 \\
{}& =
R^{(s)}_{356}(\nu) R^{(s)}_{246}(\mu) 
R^{(s)}_{145}(\lambda\nu) \
B^{(2)}_3 D^{(12)}_\lambda B^{(1)}_2   
\underline{ B^{(2)}_1
 B^{(3)}_4 D^{(23)}_{\lambda\nu} B^{(2)}_5 }
 D^{(12)}_\mu  B^{(1)}_6 
 R^{(s)}_{123}(\lambda) \\
{}& =
R^{(s)}_{356}(\nu) R^{(s)}_{246}(\mu)  \
B^{(2)}_3 D^{(12)}_\lambda B^{(1)}_2   B^{(3)}_5
 D^{(23)}_{\lambda\nu} B^{(2)}_4  B^{(3)}_1 D^{(12)}_\mu  
 B^{(1)}_6 
 R^{(s)}_{145}(\lambda\nu) R^{(s)}_{123}(\lambda) \\
{}& =
R^{(s)}_{356}(\nu) R^{(s)}_{246}(\mu)  \
B^{(2)}_3   B^{(3)}_5 D^{(23)}_{\nu} D^{(12)}_\lambda  
 D^{(23)}_{\lambda} \underline{ 
 B^{(1)}_2  B^{(2)}_4   D^{(12)}_\mu  B^{(1)}_6 }
 B^{(3)}_1
 R^{(s)}_{145}(\lambda\nu) R^{(s)}_{123}(\lambda) \\
{}& = 
R^{(s)}_{356}(\nu)  \
B^{(2)}_3   B^{(3)}_5 D^{(23)}_{\nu} D^{(12)}_\lambda  
 D^{(23)}_{\lambda}   
 B^{(2)}_6  D^{(12)}_\mu B^{(1)}_4     B^{(2)}_2  
 B^{(3)}_1
 R^{(s)}_{246}(\mu)  R^{(s)}_{145}(\lambda\nu) 
 R^{(s)}_{123}(\lambda) \\
{}& = 
R^{(s)}_{356}(\nu)   \
\underline{
B^{(2)}_3   B^{(3)}_5 D^{(23)}_{\nu}   B^{(2)}_6 }
  D^{(12)}_{\lambda\mu} B^{(1)}_4   
 D^{(23)}_{\lambda}  B^{(2)}_2   B^{(3)}_1
 R^{(s)}_{246}(\mu)  R^{(s)}_{145}(\lambda\nu) 
 R^{(s)}_{123}(\lambda) \\ 
{}& = 
B^{(3)}_6 D^{(23)}_{\nu}  B^{(2)}_5    B^{(3)}_3  
  D^{(12)}_{\lambda\mu} B^{(1)}_4   
 D^{(23)}_{\lambda}  B^{(2)}_2   B^{(3)}_1 \
 R^{(s)}_{356}(\nu) R^{(s)}_{246}(\mu)  
 R^{(s)}_{145}(\lambda\nu) R^{(s)}_{123}(\lambda) \\  
{}& =  \tilde{A}_s \
 R^{(s)}_{356}(\nu) R^{(s)}_{246}(\mu)  
 R^{(s)}_{145}(\lambda\nu) R^{(s)}_{123}(\lambda) 
 =  \tilde{A}_s \ \SM'  \,.  
\end{align*}  
The underlined terms were transformed by applying relation~(\ref{RBD}).
The remaining transformations are simply permutations of
commuting terms.
\end{proof}

\begin{proof}[\bf Proof of Theorem~\ref{RS}]

Part (i). Substitute (\ref{Rgl}) into (\ref{Rtetr2})
and pull all $\SF$'s to the left transforming the
arguments of $g$'s  according to 
(\ref{Fa}) and~(\ref{Fb}). The products of $\SF$'s
can be cancelled due to (\ref{Ftetr}) and we are left
with a product of four  $g$'s on the each side,
\begin{equation}\label{Frs2}  
 g_{q,\lambda}(\SX) \,
g_{q,\lambda\mu}( \ST) \,  
g_{q,\nu}(\SZ) \, g_{q,\mu}(\SY) =
 g_{q,\nu}( \SZ') \, g_{q,\mu}( \SY') \,
  g_{q,\lambda\nu}( \ST' )  \,
  g_{q,\lambda}( \SX')   \,, 
\end{equation} 
where the arguments of $g$'s are monomials in $a_i$'s
and $b_i$'s, for instance,
\begin{equation}\label{xyab} 
\SX=a_1^{\alpha_1} a_2^{\alpha_2} a_3^{\alpha_3}
b_1^{\beta_1} b_2^{\beta_2} b_3^{\beta_3} 
b_4^{\beta_2} 
b_5^{\beta_3+\beta_1-\beta_2} 
b_6^{-\beta_1-\beta_3} , \qquad
\SY'=a_1^{\alpha_1+\alpha_2} a_3^{\alpha_1} 
a_5^{\alpha_2} a_6^{\alpha_3}
b_2^{\beta_1} b_4^{\beta_2} b_6^{\beta_3} .
\end{equation}

If conditions (\ref{xy}) hold, then we have
$\SX=\SX'$, $\SY=\SY'$, $\SZ=\SZ'$,
$\ST= \SX\,\SY$, $\ST'= \SZ\,\SX$
and, moreover, $\SX$, $\SY$, $\SZ$ satisfy
relations (\ref{XYZ}) with $q$
replaced by $Q$ given by~(\ref{Qru}).
Therefore, by Lemma~\ref{pentdil},
$g_{q,\lambda}(t)$ must be of the form~(\ref{ru}).

Part (ii). Substitute (\ref{Rgfl}) into (\ref{Rtetr2}),
pull all $\SF$'s to the left transforming the
arguments of $g$'s and $f$'s according to 
(\ref{Fa}) and~(\ref{Fb}), and then cancel
the products of~$\SF$'s. Taking into
account that the arguments of $g$ and $f$ in (\ref{Rgfl})
satisfy conditions (\ref{xy}),
we obtain the following equality
\begin{equation}\label{Frs1}
\begin{aligned}
{}&  
 g_{q,\lambda}(\SX) f_{q,\lambda}(\tilde{\SX})
 g_{q,\lambda\mu}(\SX \SY) 
 f_{q,\lambda\mu}(\tilde{\SX} \tilde{\SY}) 
 g_{q,\nu}(\SZ)
   f_{q,\nu}(\tilde{\SZ}) g_{q,\mu}(\SY) 
   f_{q,\mu}(\tilde{\SY}) \\ 
{}& =   
  g_{q,\nu}(\SZ) f_{q,\nu}(\tilde{\SZ})
  g_{q,\mu}(\SY) f_{q,\mu}(\tilde{\SY})
  g_{q,\lambda\nu}(\SZ \SX) 
  f_{q,\lambda\nu}(\tilde{\SZ}\tilde{\SX})  
  g_{q,\lambda}(\SX) 
  f_{q,\lambda}(\tilde{\SX}) \,,
\end{aligned}
\end{equation}
where  both triples $\SX,\SY,\SZ$, and
$\tilde{\SX},\tilde{\SY},\tilde{\SZ}$
satisfy relations (\ref{XYZ}) with 
$q$ replaced, respectively, by $Q$ and $\widehat{Q}$
given by~(\ref{QQru}).
Furthermore, if conditions (\ref{xy22}) hold,
then every element of the triple $\SX,\SY,\SZ$
commutes  with  every element of the triple 
$\tilde{\SX},\tilde{\SY},\tilde{\SZ}$.
Therefore, equality (\ref{Frs1}) holds 
if its $g$ and $f$ parts satisfy separately
relation (\ref{tetr}). By Lemma~\ref{pentdil},
this requires $g_{q,\lambda}(t)$ and 
$f_{q,\lambda}(t)$ to be of the form~(\ref{ruru}).
\end{proof}

{\bf Acknowledgements.} 
We thank R. Kashaev and V. Tarasov for useful remarks.
This work was supported 
by the European Research Council (ERC) grant MODFLAT,
by the Swiss National Science Foundation grant 200020-141329,
 and in part by the Russian Fund for Basic Research 
grants  14-01-00341
and 13-01-12405-ofi-m.

{\sc \small
\noindent
Section of Mathematics, University of Geneva, 
C.P. 64, 1211 Gen\`eve 4, Switzerland \\
Steklov Mathematical Institute,
Russian Academy of Sciences,
Fontanka 27, 191023, St. Petersburg, Russia}

\end{document}